\documentclass{article}%
\usepackage{amsmath}
\usepackage{amsfonts}
\usepackage{amssymb}
\usepackage{graphicx}%
\providecommand{\U}[1]{\protect\rule{.1in}{.1in}}
\newtheorem{theorem}{Theorem}

\newtheorem{corollary}[theorem]{Corollary}

\newtheorem{definition}[theorem]{Definition}
\newtheorem{example}[theorem]{Example}

\newtheorem{lemma}[theorem]{Lemma}

\newtheorem{remark}[theorem]{Remark}

\begin{document}

\begin{center}
{\Large Quenched invariance principles for orthomartingale-like sequences}

\bigskip

Magda Peligrad and Dalibor Voln\'{y}
\end{center}

Department of Mathematical Sciences, University of Cincinnati, PO Box 210025,
Cincinnati, Oh 45221-0025, USA. \texttt{ }

Email: peligrm@ucmail.uc.edu

LMRS, CNRS-Universit\'{e} de Rouen Normandie

e-mail : dalibor.volny@univ-rouen.fr

\bigskip


\bigskip

\begin{center}
Abstract
\end{center}

In this paper we study the central limit theorem and its functional form for
random fields which are not started from their equilibrium, but rather under
the measure conditioned by the past sigma field. The initial class considered
is that of orthomartingales and then the result is extended to a more general
class of random fields by approximating them, in some sense, with an
orthomartingale. We construct an example which shows that there are
orthomartingales which satisfy the CLT but not its quenched form. This example
also clarifies the optimality of the moment conditions used for the validity
of our results. Finally, by using the so called orthomartingale-coboundary
decomposition, we apply our results to linear and nonlinear random fields.

\section{Introduction and the quenched CLT}

A very interesting type of convergence, with many practical applications, is
the almost sure conditional central limit theorem and its functional form.
This means that these theorems hold when the process is not started from its
equilibrium but it is rather started from a fixed past trajectory. In the
Markovian setting such a behavior is called a limit theorem started at a
point. In general these results are known under the name of quenched limit
theorems, as opposed to the annealed ones. A quenched CLT, for instance, is a
stronger form of convergence in distribution and implies the usual CLT. There
are examples in the literature showing that the annealed CLT\ does not
necessarily implies the quenched one. See for instance Ouchti and Voln\'{y}
(2008) and Voln\'{y} and Woodroofe (2010).

The limit theorems started at a point or from a fixed past trajectory are
often encountered in evolutions in random media and they are of considerable
importance in statistical mechanics. They are also useful for analyzing Markov
chain Monte Carlo algorithms.

In the context of random processes, this remarkable property is known for a
martingale which is stationary and ergodic, as shown in Ch. 4 in Borodin and
Ibragimov (1994) or on page 520 in Derriennic and Lin (2001). By using
martingale approximations, this result was extended to larger classes of
random variables by Cuny and Peligrad (2012), Voln\'{y} and Woodroofe (2014),
Cuny and Merlev\`{e}de (2014), among others (for a survey see Peligrad, 2015).

\quad A random field consists of multi-indexed random variables $(X_{u})_{u\in
Z^{d}}$. An important class of random fields are orthomartingales which have
been introduced by Cairoli (1969) and further developed in Khoshnevisan
(1982). They have resurfaced in many recent works. New versions of the central
limit theorem for stationary orthomartingales can be found in Wang and
Woodroofe (2013), Voln\'{y} (2015, 2019), which complement the results in Basu
and Dorea (1979), where a different definition of multiparameter martingale
was used.

In order to exploit the richness of the martingale techniques several authors
provided interesting sufficient conditions for orthomartingale approximations,
such as Gordin (2009), Voln\'{y} and Wang (2014), Cuny et al. (2015), El
Machkouri and Giraudo (2016), Peligrad and Zhang (2018 a), Giraudo (2018),
Voln\'{y} (2018). Other recent results involve random fields which are
functions of independent random variables as in El Machkouri et al. (2013) and
Wang and Woodroofe (2013). Peligrad and Zhang (2018 b) obtained necessary and
sufficient conditions for an orthomartingale approximation in the mean square.
These approximations make possible to obtain the central limit theorem
(CLT)\ for a large class of random fields. As in the case of a stochastic
processes, a natural and important question is to get a quenched version of
these CLT's. Motivated by this question, we obtain first a quenched CLT for
orthomartingales. We show by examples that the situation is different for
random fields. An orthomartingale which satisfies the CLT may fail to satisfy
the quenched CLT. The example we constructed also throws light on the
optimality of the moment conditions we use in our main result. Finally, we
extend the quenched CLT to its functional form and to a larger class of random
fields which can be decomposed into a orthomartingale and a coboundary. We
shall apply our results to linear and nonlinear random fields, often
encounters in economics.

For the sake of clarity, due to the complicated notation, we shall explain in
detail the case $d=2$ and the proof of the quenched CLT. Then, in the
subsequent sections, we shall discuss the general index set $Z^{d}$ and other
extensions of these results.

Let $(\Omega,\mathcal{K},P)$ be a probability space, let $T$ and $S$ be two
commuting, invertible, bimeasurable, measure preserving transformations from
$\Omega$ to $\Omega,$ and let $\mathcal{F}_{0,0}$ be a sub-sigma field of
$\mathcal{K}$. For all $(i,j)\in Z^{2}$ define
\begin{equation}
\mathcal{F}_{i,j}=T^{-i}S^{-j}(\mathcal{F}_{0,0})\text{, }i,j\in Z.
\label{def F}%
\end{equation}
Assume the filtration is increasing in $i$ for every $j$ fixed and increasing
in $j$ for every $i$ fixed (i.e. $\mathcal{F}_{0,0}\mathcal{\subset F}_{0,1}$
and $\mathcal{F}_{0,0}\mathcal{\subset F}_{1,0})$. For all $i$ and $j$ we also
define the following sigma algebras generated by the unions of sigma algebras:
$\mathcal{F}_{i,\infty}=\mathcal{\vee}_{m\in Z}\mathcal{F}_{i,m},$
$\mathcal{F}_{\infty,j}=\mathcal{\vee}_{n\in Z}\mathcal{F}_{n,j}$ and
$\mathcal{F}_{\infty,\infty}=\mathcal{\vee}_{n,m\in Z}\mathcal{F}_{n,m}.$ In
addition assume the filtration is commuting, in the sense that for any
integrable variable $X,$ with notation $E_{a,b}X=E(X|\mathcal{F}_{a,b}),$ we
have%
\begin{equation}
E_{u,v}E_{a,b}X=E_{a\wedge u,b\wedge v}X.\mathbb{\ } \label{com}%
\end{equation}
We introduce the stationary sequence as following. Define a function
$X_{0,0}:\Omega\rightarrow R,$ which is $\mathcal{F}_{0,0}-$measurable, and
the random field%

\begin{equation}
X_{i,j}(\omega)=X_{0,0}(T^{i}S^{j}(\omega)). \label{def X}%
\end{equation}
For the filtration $(\mathcal{F}_{i,j})$ defined by (\ref{def F})\ we call the
random field $(X_{i,j})_{i,j\in Z}$ defined by (\ref{def X})\ orthomartingale
difference field, if
\begin{equation}
E(X_{i,j}|\mathcal{F}_{u,v})=0\text{ if either }u<i\text{ or }v<j.
\label{ortho}%
\end{equation}
This definition implies that for any $i$ fixed $(X_{i,j})_{j\in Z}$ is a
sequence of martingale differences with respect to the filtration
$(\mathcal{F}_{\infty,j})_{j\in Z}$ and also for any $j$ fixed $(X_{i,j}%
)_{i\in Z}$ is a sequence of martingale differences with respect to the
filtration $(\mathcal{F}_{i,\infty})_{i\in Z}.$ Set%
\[
S_{n,v}=%
{\displaystyle\sum\nolimits_{i=0}^{n-1}}
{\displaystyle\sum\nolimits_{j=0}^{v-1}}
X_{i,j}.
\]

Below, $\Rightarrow$ denotes convergence in distribution.

The results in this paper are motivated by the following annealed CLT in
Voln\'{y} (2015), which was extended to a functional CLT in Cuny et al. (2015).

\bigskip

\textbf{Theorem A }\textit{Assume that }$(X_{i,j})_{i,j\in Z}$\textit{ is
defined by (\ref{def X}) and satisfies (\ref{ortho}). Also assume that the
filtration }$(\mathcal{F}_{i,j})_{i,j\in Z}$\textit{ is defined by
(\ref{def F}) and satisfies (\ref{com}). Assume that }$S$\textit{ (or }%
$T$\textit{) is ergodic and }$X_{0,0}$\textit{ is square integrable,
}$E(X_{0,0}^{2})=\sigma^{2}$\textit{. Then,}%
\[
\frac{1}{(nv)^{1/2}}S_{n,v}\Rightarrow N(0,\sigma^{2})\text{ \textit{when}
}n\wedge v\rightarrow\infty.
\]

Let us point out that if $S$\textit{ (or }$T$\textit{) }is ergodic,
then\textit{ }the $Z^{2}$ action generated by $S$ and $T$ is necessarily
ergodic. However the ergodicity is not enough for Theorem A to hold. In
Example 5.6 in Wang and Woodroofe (2013) and then in more detail by Voln\'{y}
(2015), a simple example of ergodic random field which does not satisfy the
central limit theorem is analyzed. Starting with two sequences of i.i.d.
random variables, centered with finite second moments, $(X_{n})$ and $(Y_{n}%
)$, the example is provided by the random field $(Z_{i,j}),$ with
$Z_{i,j}=X_{i}Y_{j}$ for all $(i,j)$.

It should be noted that Theorem A has a different area of applications than
Theorem 1 in Basu and Dorea (1979). In this latter paper the filtration is not
supposed to be commuting. For a random field $(X_{i,j})_{i,j\geq1}$ their
filtration $(\mathcal{K}_{n,m})$ is generated by the variables $\{X_{i,j}%
:(j\geq1,1\leq i\leq n)\cup(i\geq1,1\leq j\leq m)\}.$ Suppose $(\xi_{i,j})$
are i.i.d., standard normal random variables. Then, Theorem A can be applied,
for instance, to the random field $(X_{i,j}),$ where $X_{i,j}(\omega
)=X_{0,0}(T^{i}S^{j}(\omega))$ with $X_{0,0}=\xi_{-1,0}\xi_{0,-1}$ and
$\mathcal{F}_{0,0}=\sigma(\xi_{i,j},i\leq0,j\leq0)$ but the result in Basu and
Dorea (1979) cannot. On the other hand the random field $(Y_{i,j}),$ defined
by $Y_{i,j}=Y_{0,0}(T^{i}S^{j}(\omega))$ with $Y_{0,0}=\sum\nolimits_{k=1}%
^{\infty}a_{k}(\xi_{k,0}+\xi_{0,k})$ and $\sum\nolimits_{k=1}^{\infty}%
|a_{k}|<\infty,$ can be treated by the result in Basu and Dorea (1979) but not
by Theorem A.

It should also be noted that Theorem A allows to study the central limit
theorem for orthomartingales which are not defined by a Bernoulli $Z^{2}$-action.

The aim of this paper is to establish a quenched version of Theorem A.

We denote by $P^{\omega}(\cdot)=P_{0,0}^{\omega}(\cdot)$ a version of the
regular conditional probability $P(\cdot|\mathcal{F}_{0,0})(\omega)$.

One of the results of this paper is the following theorem:

\begin{theorem}
\label{QCLTmart}Assume that $(X_{i,j})_{i,j\in Z}$ is defined by (\ref{def X})
and satisfies (\ref{ortho}). Also assume that the filtration $(\mathcal{F}%
_{i,j})_{i,j\in Z}$ is defined by (\ref{def F}) and satisfies (\ref{com}).
Assume that $S$ (or $T$) is ergodic and $X_{0,0}$ is square integrable,
$E(X_{0,0}^{2})=\sigma^{2}$. Then, for $P$-almost all $\omega\in\Omega,$%
\begin{equation}
\frac{1}{n}S_{n,n}\Rightarrow N(0,\sigma^{2})\text{ under }P^{\omega}.
\label{QCLTn}%
\end{equation}
In addition, if
\begin{equation}
E(X_{0,0}^{2}\log(1+|X_{0,0}|))<\infty, \label{moment condition}%
\end{equation}
then for almost all all $\omega\in\Omega,$
\begin{equation}
\frac{1}{(nv)^{1/2}}S_{n,v}\Rightarrow N(0,\sigma^{2})\text{ under }P^{\omega
}\text{ when }n\wedge v\rightarrow\infty. \label{QCLTnm}%
\end{equation}

\end{theorem}

We would like to mention that, because by integration the quenched CLT implies
the annealed CLT, the conclusion in Theorem \ref{QCLTmart} implies the CLT in
Theorem A. However, when the summation on the rectangles is not restricted,
the integrability assumption (\ref{moment condition})\ is stronger than in
Theorem A. Later on, in Theorem \ref{Th funct CLT}, we shall extend this
result to a functional central limit theorem. \ Let us also notice that the
second part of Theorem \ref{QCLTmart} does not always hold under the
assumption $E(X_{0,0}^{2})<\infty$. As a matter of fact we are going to
provide an example to support this claim.

\begin{theorem}
\label{Counterex}Under the setting used in Theorem \ref{QCLTmart}, there is a
stationary sequence $(X_{n,m})_{n,m\in Z}$\ satisfying (\ref{ortho}), adapted
to a commuting filtration $(\mathcal{F}_{i,j})_{i,j\in Z},$ with
$E(X_{0,0}^{2}\ln(1+|X_{0,0}|))=\infty,$ for any $0<\varepsilon<1,$
$E(X_{0,0}^{2}\ln^{1-\varepsilon}(1+|X_{0,0}|))<\infty$ and such that
$(S_{n,m}/\sqrt{nm})_{(n,m)\in Z^{2}}$ does not satisfy the quenched CLT in
(\ref{QCLTnm}).
\end{theorem}

We mention that, as a matter of fact, in our examples, both transformations
constructed for the definition of $(X_{n,m})_{n,m\in Z}$ and for the
filtration $(\mathcal{F}_{i,j})_{i,j\in Z},$ are ergodic. Also, this example
satisfies the quenched CLT in (\ref{QCLTn}).

The detailed proofs of these two theorems are contained in Section
\ref{proof main}. Various extensions of Theorem \ref{QCLTmart} will be given
in subsequent sections.

In Section \ref{functional CLT} we formulate the functional form of the
quenched CLT and we indicate how to prove it, by adapting the arguments from
the proof of Theorem \ref{QCLTmart} and some other proofs of several known results.

For the sake of applications, in Section \ref{cob}, we extend the results
beyond orthomartingales, to a class of random fields which can be decomposed
into an orthomartingale and a generalized coboundary.

In Section \ref{section gend} we show that Theorem \ref{QCLTmart} remains
valid for random fields indexed by $Z^{d},$ $d>2.$ The only difference is that
we replace condition (\ref{moment condition}) by $E(X_{0,0}^{2}\log
^{d-1}(1+|X_{0,0}|))<\infty.$

In Section \ref{Appl} we apply our results to linear and nonlinear random
fields with independent innovations. Several useful results for our proofs are
given in Section 7.

\section{Proofs of Theorems \ref{QCLTmart} and \ref{Counterex}
\label{proof main}}

\textbf{Proof of Theorem \ref{QCLTmart}}

\bigskip

To fix the ideas, let us suppose that the transformation $S$ is ergodic. Let
us denote by $\hat{T}$ and $\hat{S}$ the operators on $L_{2}$, defined by
$\hat{T}f=f\circ T$ and $\hat{S}f=f\circ S.$ Everywhere in the paper, for $x$
real, we shall denote by $[x]$ the integer part of $x.$

By using a truncation argument, we show first that, without restricting the
generality, we can prove the theorem under the additional assumption that the
variables are bounded. We shall introduce the following projection operators:%
\[
\mathcal{P}_{i,j}(X)=E_{i,j}(X)-E_{i,j-1}(X)-E_{i-1,j}(X)+E_{i-1,j-1}(X).
\]
Let $A$ be a positive integer. Denote $X_{i,j}^{^{\prime}}=X_{i,j}%
I(|X_{i,j}|\leq A)$ and $X_{i,j}^{"}=X_{i,j}I(|X_{i,j}|>A).$ Therefore, we can
represent $(X_{i,j})$ as a sum of two orthomartingale differences adapted to
the same filtration.
\begin{equation}
X_{i,j}=\mathcal{P}_{i,j}(X_{i,j}^{^{\prime}})+\mathcal{P}_{i,j}(X_{i,j}^{"}).
\label{truncation}%
\end{equation}
Note that,%
\[
|\mathcal{P}_{0,0}(X_{0,0}^{"})|\leq|X_{0,0}|+E_{-1,0}|X_{0,0}|+E_{0,-1}%
|X_{0,0}|+E_{-1,1}|X_{0,0}|.
\]
Whence, by the properties of conditional expectation, $E(X_{0,0})^{2}<\infty$
implies
\begin{equation}
E(\mathcal{P}_{0,0}(X_{0,0}^{"}))^{2}<\infty\label{secmp}%
\end{equation}
and $E(X_{0,0}^{2}\log(1+|X_{0,0}|))<\infty$ implies
\begin{equation}
E((\mathcal{P}_{0,0}(X_{0,0}^{"}))^{2}\log(1+|(\mathcal{P}_{0,0}(X_{0,0}%
^{"}))|)<\infty. \label{secmp2}%
\end{equation}
Set%
\[
S_{n,v}^{^{\prime}}=%
{\displaystyle\sum\nolimits_{i=0}^{n-1}}
{\displaystyle\sum\nolimits_{j=0}^{v-1}}
\mathcal{P}_{i,j}(X_{i,j}^{^{\prime}})\text{ and }S_{n,v}^{"}=%
{\displaystyle\sum\nolimits_{i=0}^{n-1}}
{\displaystyle\sum\nolimits_{j=0}^{v-1}}
\mathcal{P}_{i,j}(X_{i,j}^{"}).
\]
We shall show that, for $P-$almost all $\omega,$
\[
\lim_{A\rightarrow\infty}\text{ }\lim\sup_{n\wedge v\rightarrow\infty
}P^{\omega}(\frac{1}{(nv)^{1/2}}|S_{n,v}^{"}|>\varepsilon)=0.
\]
By conditional Markov inequality, it is enough to show that
\begin{equation}
\lim_{A\rightarrow\infty}\lim_{n\wedge v\rightarrow\infty}\frac{1}{nv}%
E_{0,0}(S_{n,v}^{"})^{2}=0\text{ \ \ a.s.} \label{negli1}%
\end{equation}
By the orthogonality of the orthomartingale differences, we have that
\begin{equation}
\frac{1}{nv}E_{0,0}((S_{n,v}^{"})^{2})=\frac{1}{nv}%
{\displaystyle\sum\nolimits_{i=0}^{n-1}}
{\displaystyle\sum\nolimits_{j=0}^{v-1}}
E_{0,0}(\mathcal{P}_{i,j}(X_{i,j}^{"}))^{2}. \label{orthogonality}%
\end{equation}
Note that the conditional expectation introduces a family of operators defined
by
\[
Q_{1}(f)=E_{0,\infty}(\hat{T}f)\text{ ; }Q_{2}(f)=E_{\infty,0}(\hat{S}f).
\]
So, using (\ref{com}), we can write
\[
E_{0,0}(\mathcal{P}_{i,j}(X_{i,j}^{"}))^{2}=Q_{1}^{i}Q_{2}^{j}(\mathcal{P}%
_{0,0}(X_{0,0}^{"}))^{2}.
\]

Since $Q_{1}$ and $Q_{2}$ are integral preserving Dunford-Schwartz operators,
by the ergodic theorem (see Theorem 3.5 in Ch. 6 in Krengel, 1985), if we
assume finite second moment, by (\ref{secmp}),%
\[
\lim_{n\rightarrow\infty}\frac{1}{n^{2}}%
{\displaystyle\sum\nolimits_{i=0}^{n-1}}
{\displaystyle\sum\nolimits_{j=0}^{n-1}}
Q_{1}^{i}Q_{2}^{j}(\mathcal{P}_{0,0}(X_{0,0}^{"}))^{2}=E(\mathcal{P}%
_{0,0}(X_{0,0}^{"}))^{2}\text{ a.s.}%
\]
If we assume $E(X_{0,0}^{2}\log(1+|X_{0,0}|))<\infty$ then, by (\ref{secmp2})
and Theorem 1.1 in Ch. 6, Krengel (1985), we obtain
\begin{equation}
\lim_{n\wedge v\rightarrow\infty}\frac{1}{nv}%
{\displaystyle\sum\nolimits_{i=0}^{n-1}}
{\displaystyle\sum\nolimits_{j=0}^{v-1}}
Q_{1}^{i}Q_{2}^{j}(\mathcal{P}_{0,0}(X_{0,0}^{"}))^{2}=E(\mathcal{P}%
_{0,0}(X_{0,0}^{"}))^{2}\text{ a.s.} \label{erg}%
\end{equation}
Clearly $\lim_{A\rightarrow\infty}\mathcal{P}_{0,0}(X_{0,0}^{"})=0$ a.s. So,
by the dominated convergence theorem,
\[
\lim_{A\rightarrow\infty}E(\mathcal{P}_{0,0}(X_{0,0}^{"}))^{2}=0,
\]
and (\ref{negli1}) is established. By Theorem 3.2 in Billingsley (1999), in
order to establish conclusion (\ref{QCLTnm}) of Theorem \ref{QCLTmart}, it is
enough to show that for $A$ fixed, for almost all $\omega\in\Omega,$%
\[
\frac{1}{(nv)^{1/2}}S_{n,v}^{\prime}\Rightarrow N(0,\sigma_{A}^{2})\text{
under }P^{\omega}\text{ as }n\wedge v\rightarrow\infty,\text{ and }\sigma
_{A}^{2}\rightarrow\sigma^{2}\text{ as }A\rightarrow\infty.
\]
Above, $\sigma_{A}^{2}=E(\mathcal{P}_{0,0}(X_{0,0}^{^{\prime}}))^{2}.$
Clearly, when $A\rightarrow\infty,$ $\sigma_{A}^{2}\rightarrow\sigma^{2}.$
Therefore the result is established if we prove Theorem \ref{QCLTmart}\textbf{
}for orthomartingale differences which are additionally uniformly bounded.

So, in the rest of the proof, without restricting the generality, we shall
assume that the variables $(X_{i,j})_{i,j\in Z}$ are bounded by a positive
constant $C$. Also, proving the result for $n>v\rightarrow\infty$ is
equivalent to proving it for any subsequence $(n,v_{n})$ with $v_{n}%
\rightarrow\infty$ as $n\rightarrow\infty$. To ease the notation we shall
denote $v=v_{n}$.

Denote
\begin{equation}
F_{i,v}=\frac{1}{v^{1/2}}%
{\displaystyle\sum\nolimits_{j=0}^{v-1}}
X_{i,j}. \label{defFiv}%
\end{equation}
We treat the double summation as a sum of a triangular array of martingale
differences $(F_{i,v})_{i\geq0}:$%
\[
\frac{1}{(nv)^{1/2}}S_{n,v}=\frac{1}{n^{1/2}}%
{\displaystyle\sum\nolimits_{i=0}^{n-1}}
F_{i,v}.
\]
We shall apply Theorem 1 in G\"{a}nssler and H\"{a}usler (1979), given for
convenience in Theorem\ \ref{cltGH} from Section \ref{AUX}, to $D_{n,i}%
=F_{i,v}/\sqrt{n}$. We have to show that for almost all $\omega,$ both
conditions of this theorem are satisfied. Namely we shall verify that $P-$for
almost all $\omega\in\Omega$ and all rationals $q\in\lbrack0,1]$
\begin{equation}
\lim_{n\rightarrow\infty}\frac{1}{n}E_{0,0}|\sum\nolimits_{i=0}^{[(n-1)q]}%
(F_{i,v}^{2}-\sigma^{2})|=0\text{. } \label{limit cond}%
\end{equation}
and%
\begin{equation}
\frac{1}{n}E_{0,0}\max_{0\leq i\leq n-1}F_{i,v}^{2}\text{ is bounded.}
\label{max bound cond}%
\end{equation}

We verify first (\ref{limit cond}). Note that, since the rationals are
countable, it is enough to show that for any $q$ rational%
\[
\lim_{n\rightarrow\infty}\frac{1}{n}E_{0,0}|\sum\nolimits_{i=0}^{[(n-1)q]}%
(F_{i,v}^{2}-\sigma^{2})|=0\text{ \ }P-\text{a.s.}%
\]
\ We verify it first with $q=1$ and use a blocking procedure.

Let $m\geq1$ be a fixed integer and define consecutive blocks of indexes of
size $m$, $I_{j}(m)=\{(j-1)m,...,mj-1\}.$ In the set of integers from $0$ to
$n-1$ we have $u=u_{n}(m)=[n/m]$ such blocks of integers and a last one
containing less than $m$ indexes. Practically, by the triangle inequality, we
write%
\begin{gather*}
\frac{1}{n}|\sum\nolimits_{i=0}^{n-1}(F_{i,v}^{2}-\sigma^{2})|\leq\\
\frac{1}{n}\sum\nolimits_{j=1}^{u}|\sum\nolimits_{k\in I_{j}(m)}(F_{k,v}%
^{2}-\sigma^{2})|+\frac{1}{n}|\sum\nolimits_{k=um}^{n-1}(F_{k,v}^{2}%
-\sigma^{2})|\leq\\
\frac{1}{u}\sum\nolimits_{j=1}^{u}|\frac{1}{m}\sum\nolimits_{k\in I_{j}%
(m)}F_{k,v}^{2}-\sigma^{2}|+\frac{1}{n}|\sum\nolimits_{k=um}^{n-1}(F_{k,v}%
^{2}-\sigma^{2})|=\\
I_{n,m}+II_{n,m}.
\end{gather*}
The task is now to show that
\begin{equation}
\lim_{m\rightarrow\infty}\lim_{n\wedge v\rightarrow\infty}E_{0,0}%
(I_{n,m})=0\text{ }\ \text{a.s.} \label{one}%
\end{equation}
and%
\begin{equation}
\lim_{m\rightarrow\infty}\lim_{n\wedge v\rightarrow\infty}E_{0,0}%
(II_{n,m})=0\text{ }\ \text{a.s.} \label{two}%
\end{equation}
Let us treat first the limit of $E_{0,0}(I_{n,m})$. Let $N_{0}$ be a fixed
integer and consider $n\wedge v>N_{0}$. By using the properties of the
conditional expectations and (\ref{com}) we obtain the following bound for
$E_{0,0}(I_{n,m}):$
\begin{gather*}
E_{0,0}(I_{n,m})=\frac{1}{u}E_{0,0}\sum\nolimits_{j=1}^{u}|\frac{1}{m}%
\sum\nolimits_{k\in I_{j}(m)}F_{k,v}^{2}-\sigma^{2}|\\
=\frac{1}{u}E_{0,0}\sum\nolimits_{j=1}^{u}E_{(j-1)m,0}|\frac{1}{m}%
\sum\nolimits_{k\in I_{j}(m)}F_{k,v}^{2}-\sigma^{2}|\\
=E_{0,0}\frac{1}{u}\sum\nolimits_{i=0}^{u-1}\hat{T}^{im}E_{0,0}|\frac{1}%
{m}\sum\nolimits_{k=0}^{m-1}F_{k,v}^{2}-\sigma^{2}|\\
\leq E_{0,0}\frac{1}{u}\sum\nolimits_{i=0}^{u-1}\hat{T}^{im}(h_{m,N_{0}}),
\end{gather*}
where we have used the notation%
\[
h_{m,N_{0}}=\sup_{v>N_{0}}E_{0,0}|\frac{1}{m}\sum\nolimits_{k=0}^{m-1}%
F_{k,v}^{2}-\sigma^{2}|.
\]
Note that $h_{m,N_{0}}$ is bounded. Indeed, by the martingale property and the
uniform boundedness of the variables by $C$, it follows that
\begin{gather*}
h_{m,N_{0}}\leq\sigma^{2}+\frac{1}{m}\sum\nolimits_{k=0}^{m-1}\sup_{v>N_{0}%
}E_{0,0}(F_{k,v}^{2})\\
=\sigma^{2}+\frac{1}{m}\sum\nolimits_{k=0}^{m-1}\sup_{v>N_{0}}E_{0,0}(\frac
{1}{v}\sum\nolimits_{u=0}^{v-1}X_{k,u}^{2})\leq\sigma^{2}+C^{2}.
\end{gather*}
By the ergodic theorem, (see Theorem 11.4 in Eisner et al., 2015 or Corollary
3.8 in Ch. 3, Krengel, 1985) for each $m$ and $N_{0}$
\[
\lim_{u\rightarrow\infty}\frac{1}{u}\sum\nolimits_{i=0}^{u-1}\hat{T}%
^{im}h_{m,N_{0}}=E(h_{m,_{N_{0}}}|I)=E_{I}(h_{m,_{N_{0}}})\text{ a.s.,}%
\]
where $I$ is the invariant sigma field for the operator $T$. Furthermore, we
also have that
\[
\frac{1}{u}\sum\nolimits_{i=0}^{u-1}\hat{T}^{im}h_{m,N_{0}}\leq\sigma
^{2}+C^{2}.
\]
So, by Theorem 34.2 (v) in Billingsley (1995) (see Theorem \ref{Bill} in
Section \ref{AUX}) we derive that
\[
\lim_{u\rightarrow\infty}E_{0,0}\frac{1}{u}\sum\nolimits_{i=0}^{u-1}\hat
{T}^{im}h_{m,N_{0}}=E_{0,0}E_{I}(h_{m,_{N_{0}}})\text{ a.s.}%
\]
Since the functions are bounded, by applying twice, consecutively, Theorem
\ref{Bill}, we obtain that%
\[
\lim_{N_{0}\rightarrow\infty}\lim_{u\rightarrow\infty}E_{0,0}\frac{1}{u}%
\sum\nolimits_{i=0}^{u-1}\hat{T}^{im}h_{m,N_{0}}=E_{0,0}E_{I}(\lim
_{N_{0}\rightarrow\infty}h_{m,_{N_{0}}})\text{ a.s.}%
\]
Clearly, because the variables are bounded, for every $m$ fixed
\begin{align*}
E_{0,0}E_{I}(\lim_{N_{0}\rightarrow\infty}h_{m,_{N_{0}}})  &  =E_{0,0}%
E_{I}(\lim\sup_{v}E_{0,0}|\frac{1}{m}\sum\nolimits_{k=0}^{m-1}F_{k,v}%
^{2}-\sigma^{2}|)\\
&  \leq E_{0,0}E_{I}E_{0,0}(\lim\sup_{v}E_{\infty,0}|\frac{1}{m}%
\sum\nolimits_{k=0}^{m-1}F_{k,v}^{2}-\sigma^{2}|).
\end{align*}
Now, by using again the fact that the variables are bounded and using Theorem
\ref{Bill}, in order to show that%
\[
\lim_{m\rightarrow\infty}E_{0,0}E_{I}(\lim_{N_{0}\rightarrow\infty
}h_{m,_{N_{0}}})=0\text{ }P\text{-a.s.}%
\]
it is enough to show that
\begin{equation}
\lim_{m\rightarrow\infty}\lim\sup_{v}E_{\infty,0}|\frac{1}{m}\sum
\nolimits_{k=0}^{m-1}F_{k,v}^{2}-\sigma^{2}|=0\text{ }\ \text{a.s.}
\label{to show}%
\end{equation}
With this aim, we note first that by the ergodicity of $S$ and the fact that
the variables are bounded, it follows that, for any $k,$
\begin{equation}
\lim_{v\rightarrow\infty}E_{\infty,0}F_{k,v}^{2}=\lim_{v\rightarrow\infty
}\frac{1}{v}E_{\infty,0}(\sum\nolimits_{j=0}^{v-1}X_{k,j}^{2})=\sigma^{2}.
\label{erg2}%
\end{equation}
Denote $P_{\infty,0}^{\omega}(\cdot)=P(\cdot|\mathcal{F}_{\infty,0}).$ We also
know that for any $k$, by the quenched CLT for stationary martingale
differences (see, for instance, Ch. 4 in Borodin and Ibragimov (1994) or
Derrienic and Lin (2001)), for almost all $\omega,$ $F_{k,v}\Rightarrow N_{k}$
under $P_{\infty,0}^{\omega}$, where $N_{k}$ is a centered normal random
variable with variance $\sigma^{2}.$ Therefore, by the sufficiency part of the
convergence of moments associated to weak convergence, namely Theorem 3.6\ in
Billingsley (1999), we have that%
\begin{equation}
(F_{k,v}^{2})_{v\geq1}\text{ is uniformly integrable under}\ P_{\infty
,0}^{\omega}\text{ for almost all }\omega. \label{U.I}%
\end{equation}
By the functional quenched CLT\ for martingales (see Ch. 4 in Borodin and
Ibragimov (1994)), for almost all $\omega$,$\ $we know that
\[
(F_{0,v},F_{1,v},...,F_{m-1,v})\Rightarrow(N_{0},N_{1},...,N_{m-1})\text{
under }P_{\infty,0}^{\omega}\text{ as }v\rightarrow\infty,
\]
where $(N_{0},N_{1},...,N_{m-1})$ is a Gaussian vector of centered normal
variables with variance $\sigma^{2}$. But since $(F_{j,v})_{j\in Z}$ are
uncorrelated it follows by (\ref{U.I}) that the variables in $(N_{i})_{i\geq
0}$ are also uncorrelated and therefore $(N_{i})_{i\geq0}$ is an i.i.d.
sequence. By the continuous mapping theorem,
\[
\frac{1}{m}\sum\nolimits_{k=0}^{m-1}(F_{k,v}^{2}-\sigma^{2})\Rightarrow
\frac{1}{m}\sum\nolimits_{k=0}^{m-1}(N_{k}^{2}-\sigma^{2})\text{ under
}P_{\infty,0}^{\omega}\text{ for almost all }\omega\text{.}%
\]
By (\ref{U.I}) it follows that $(\sum\nolimits_{k=0}^{m-1}(F_{k,v}^{2}%
-\sigma^{2}))_{v\geq1}$ is also uniformly integrable, so we can apply the
convergence of moments from Theorem 3.5\ in Billingsley (1999). Therefore,
denoting by $\mathcal{E}$ the expectation in rapport with the probability on
the space where the variables $(N_{k})^{\prime}s$ are defined, we obtain
\[
\lim_{v\rightarrow\infty}E_{\infty,0}|\frac{1}{m}\sum\nolimits_{k=0}%
^{m-1}(F_{k,v}^{2}-\sigma^{2})|=\mathcal{E}|\frac{1}{m}\sum\nolimits_{k=0}%
^{m-1}(N_{k}^{2}-\sigma^{2})|\text{ \ }\ \text{a.s.}%
\]
By letting $m\rightarrow\infty$ and using the law of large numbers for an
i.i.d. sequence, we obtain%
\[
\lim_{m\rightarrow\infty}\mathcal{E}(|\frac{1}{m}\sum\nolimits_{k=0}%
^{m-1}(N_{k}^{2}-\sigma^{2})|=0.
\]
Therefore (\ref{to show}) follows. As a consequence, we obtain (\ref{one}).

In order to treat the term (\ref{two}), we estimate%
\begin{gather*}
E_{0,0}(II_{n,m})=E_{0,0}\frac{1}{n}|\sum\nolimits_{k=um}^{n-1}(F_{k,v}%
^{2}-\sigma^{2})|\leq\frac{m}{n}\sigma^{2}+E_{0,0}\frac{1}{n}\sum
\nolimits_{k=um}^{n-1}F_{k,v}^{2}\\
\leq\frac{m}{n}\sigma^{2}+\frac{1}{n}\sum\nolimits_{k=um}^{n-1}\frac{1}{v}%
{\displaystyle\sum\nolimits_{j=0}^{v-1}}
E_{0,0}X_{k,j}^{2}\leq\frac{m}{n}(\sigma^{2}+C^{2})\text{ a.s.}%
\end{gather*}
Whence, (\ref{two})\ follows, by passing to the limit first with
$n\rightarrow\infty$ followed by $m\rightarrow\infty$.

Overall, we have shown that%
\[
\lim_{n\wedge v\rightarrow\infty}\frac{1}{n}E_{0,0}|\sum\nolimits_{u=0}%
^{n-1}(F_{u,v}^{2}-\sigma^{2})|=0\text{ a.s.}%
\]
If we replace now $n-1$ by $[(n-1)q]$, with $q$ a rational number, we easily
see that we also have convergence to $q\sigma^{2}$ and (\ref{limit cond}) follows.

It remains to verify the second condition of Theorem \ref{cltGH}, namely to
prove (\ref{max bound cond}). To show it, note that, by the martingale
property,
\begin{align*}
\frac{1}{n}E_{0,0}(\max_{0\leq i\leq n-1}F_{i,v}^{2})  &  \leq\frac{1}%
{n}E_{0,0}(%
{\displaystyle\sum\nolimits_{i=0}^{n-1}}
F_{i,v}^{2})\\
&  =\frac{1}{nv}(%
{\displaystyle\sum\nolimits_{i=0}^{n-1}}
{\displaystyle\sum\nolimits_{u=0}^{v-1}}
E_{0,0}(X_{i,u}^{2}))\leq C^{2}\text{ a.s.}%
\end{align*}
The proof of the theorem is now complete. $\square$

\bigskip

\textbf{Proof of Theorem \ref{Counterex}}

\bigskip

We start with an i.i.d. random field $(\xi_{n,m})_{n,m\in Z}$ defined on a
probability space $(\Omega,\mathcal{K},P)$ with the distribution
\begin{equation}
P(\xi_{0,0}=-1)=P(\xi_{0,0}=1)=1/2. \label{defcsi}%
\end{equation}
Without restricting the generality we shall define $(\xi_{\mathbf{u}%
})_{\mathbf{u}\in Z^{2}}$ in a canonical way on the probability space $\Omega$
$=R^{Z^{2}}$, endowed with the $\sigma-$field $\mathcal{B},$ generated by
cylinders. Then, if $\omega=(x_{\mathbf{v}})_{\mathbf{v}\in Z^{2}}$, we define
$\xi_{\mathbf{u}}^{\prime}(\omega)=x_{\mathbf{u}}$. We construct a probability
measure $P^{\prime}$ on $\mathcal{B}$ such that for all $B\in\mathcal{B}$, any
$m$ and $\mathbf{u}_{1},...,\mathbf{u}_{m}$ we have%
\[
P^{\prime}((x_{\mathbf{u}_{1}},...,x_{\mathbf{u}_{m}})\in B)=P((\xi
_{\mathbf{u}_{1}},...,\xi_{\mathbf{u}_{m}})\in B).
\]
The new sequence $(\xi_{\mathbf{u}}^{\prime})_{\mathbf{u}\in Z^{2}}$ is
distributed as $(\xi_{\mathbf{u}})_{\mathbf{u}\in Z^{2}}$ and re-denoted by
$(\xi_{\mathbf{u}})_{\mathbf{u}\in Z^{2}}$. We shall also re-denote
$P^{\prime}$ as $P.$ Now on $R^{Z^{2}}$ we introduce the operators%
\[
T_{\mathbf{u}}((x_{\mathbf{v}})_{\mathbf{v}\in Z^{2}})=(x_{\mathbf{v+u}%
})_{\mathbf{v}\in Z^{2}}.
\]
Two of them will play an important role, namely when $\mathbf{u=}(1,0)$ and
when $\mathbf{u=}(0,1).$ By interpreting the indexes as notations for the
lines and columns of a matrix, we shall call%

\[
T((x_{u,v})_{(u,v)\in Z^{2}})=(x_{u+1,v})_{(u,v)\in Z^{2}}%
\]
the vertical shift and%
\[
S((x_{u,v})_{(u,v)\in Z^{2}})=(x_{u,v+1})_{(u,v)\in Z^{2}}%
\]
the horizontal shift. Introduce the filtration $\mathcal{F}_{n,m}=\sigma
(\xi_{i,j},i\leq n,j\leq m)$ and notice that this filtration is commuting. We
assume $\mathcal{K}=\mathcal{F}_{\infty,\infty}.$ The transformations $T$ and
$S$ are invertible, measure preserving, commuting and ergodic. Furthermore
$T_{i,j}=T^{i}S^{j}.$

For a measurable function $f$ defined on $R^{Z^{2}}$ define
\begin{equation}
X_{j,k}=f(T^{j}S^{k}(\mathbf{\xi}_{a,b})_{a\leq0,b\leq0}). \label{defXfield}%
\end{equation}
We notice that the variables are adapted to the filtration $(\mathcal{F}%
_{n,m})_{n,m\in Z}$.

\bigskip

As an important step for constructing our example we shall establish the
following lemma:

\begin{lemma}
\label{Rohling2}For every $n$ and every $\varepsilon>0$ we can find a set
$F=F(n,\varepsilon)$ which is $\mathcal{F}_{0,0}$ measurable and such that%
\[
P(F)\geq\frac{1}{n^{2}}(1-\varepsilon).
\]
Furthermore, for any $0\leq i,j\leq n-1,$ $0\leq k,\ell\leq n-1$ with
$(i,j)\neq(k,\ell)$ we have
\begin{equation}
P(T_{i,j}^{-1}F\cap T_{k,l}^{-1}F)=0. \label{disjoint}%
\end{equation}

\end{lemma}

\textbf{Proof of Lemma \ref{Rohling2}}.

\bigskip

Let $n$ be an integer and let $\varepsilon>0.$ By using Rokhlin lemma (see
Theorem \ref{Rohling} in Section \ref{AUX}), construct $B\in\mathcal{K}$ with
\begin{equation}
P(B)\geq(1-\frac{\varepsilon}{2})\frac{1}{n^{2}} \label{estB}%
\end{equation}
and for $0\leq i,j\leq n-1$, $T_{i,j}^{-1}B$ are disjoint for distinct pair of
indexes. Since $\mathcal{K}$ is generated by the field $\cup_{\mathbf{n}%
}\mathcal{F}_{\mathbf{n}},$ we can find a set $E$ in $\cup_{\mathbf{n}%
}\mathcal{F}_{\mathbf{n}}$ such that
\begin{equation}
P(B\Delta E)<\frac{\varepsilon}{8n^{4}}. \label{sym}%
\end{equation}
Since $E$ belongs to $\cup_{\mathbf{n}}\mathcal{F}_{\mathbf{n}},$ there is a
$\mathbf{m}$ such that $E\in\mathcal{F}_{\mathbf{m}}.$ So $T_{\mathbf{m}%
}(E)\in\mathcal{F}_{\mathbf{0}}.$ Denote $G=T_{\mathbf{m}}(E)$ and set%
\[
F=G\setminus\cup_{(i,j)\in D}T_{i,j}^{-1}G,
\]
where $D=\{0\leq i,j\leq n-1,(i,j)\neq(0,0)\}$. Note now that for all
$(i,j)\in D,$%
\[
P(F\cap T_{i,j}^{-1}F)=0,
\]
which implies (\ref{disjoint}). Also, by stationarity,
\[
P(F)=P(E)-P(E\cap(\cup_{(i,j)\in D}T_{i,j}^{-1}E))\geq P(E)-\sum
\nolimits_{(i,j)\in D}P(E\cap T_{i,j}^{-1}E).
\]
But for $(i,j)\in D,$
\[
P(E\cap T_{i,j}^{-1}E)\leq2P(E\setminus B)\leq\frac{\varepsilon}{4n^{4}}.
\]
Therefore, by the above considerations, (\ref{sym}) and (\ref{estB}) we obtain%
\[
P(F)\geq P(E)-\frac{\varepsilon}{4n^{2}}\geq P(B)-\frac{\varepsilon}{8n^{4}%
}-\frac{\varepsilon}{4n^{2}}\geq\frac{1-\varepsilon}{n^{2}}.
\]
$\square$

\bigskip

Next, we obtain a lemma which is the main step in the construction of the
example. In the sequel, we use the notation $a_{n}\sim b_{n}$ for
$\lim_{n\rightarrow\infty}a_{n}/b_{n}=1.$

\begin{lemma}
\label{NTight}There is a strictly stationary random field of integrable
positive random variables $(U_{i,j})_{i,j\in Z},$ coordinatewise ergodic, such
that for any $0<\varepsilon<1,$ $E|U_{0,0}|\ln^{1-\varepsilon}(1+|U_{0,0}%
|)<\infty$ and such that for almost all $\omega,$ $(U_{n,v}/nv)_{n,v\in Z}$ is
not tight under $P^{\omega}.$
\end{lemma}

\textbf{Proof of Lemma \ref{NTight}}.

\bigskip

By Lemma \ref{Rohling2}, for $n\geq2$ and $\varepsilon=1/2,$ we can find sets
$F_{n}\in\mathcal{F}_{-n,-n}$ such that $P(F_{n})=1/2n^{2}$ and such that for
any $0\leq i,j\leq n-1,$ $0\leq k,\ell\leq n-1$ with $(i,j)\neq(k,\ell)$ we
have $P(T_{i,j}^{-1}F_{n}\cap T_{k,l}^{-1}F_{n})=0.$

Now, we consider independent copies of the probability space $(\Omega
,\mathcal{K},P),$ denoted by $(\Omega^{(m)},\mathcal{K}^{(m)},P^{(m)}%
)_{m\geq1},$ and introduce the product space~$\mathbf{\Omega}=\prod
_{m=1}^{\infty}\Omega^{(m)}$ endowed with the sigma algebra generated by
cylinders, $\mathbf{K}=\prod_{m=1}^{\infty}\mathcal{K}^{(m)}$. We also
introduce on $\mathbf{K}$ the product probability $\mathbf{P}=\prod
_{m=1}^{\infty}P^{(m)},$ $P^{(m)}=P$. In this space consider sets $F_{n}%
^{(n)}$ which are products of $\Omega$ with the exception of the $n$-th
coordinate which is $F_{n}.$

On $\mathbf{\Omega,}$ define a random variable $f_{n}$ by the following
formula:
\begin{equation}
f_{n}=\frac{n}{\ln^{2}n}1_{F_{n}^{(n)}.} \label{def fn}%
\end{equation}
Let $A_{n}$ be the following event:
\[
A_{n}=\{\text{there are }i,j\text{, }\ln n\leq i,j\leq n-1,\text{ such that
}f_{n}\circ\mathbf{T}_{i,j}/ij\geq1\}.
\]
where $\mathbf{T}_{i,j}=(T_{i,j},T_{i,j},...).$ Since $f_{n}\circ
\mathbf{T}_{i,j}$ is $\prod\nolimits_{m=1}^{\infty}\mathcal{F}_{0,0}^{(m)}$
measurable, for $\omega\in A_{n},$ there are $i,j$, $\ln n\leq i,j\leq n-1,$
such that
\begin{equation}
\mathbf{P}^{\omega}(f_{n}\circ\mathbf{T}_{i,j}/ij\geq1)=1. \label{basic}%
\end{equation}
Note now that $f_{n}\circ\mathbf{T}_{i,j}/ij\geq1$ if and only if
$1_{F_{n}^{(n)}}\circ\mathbf{T}_{i,j}\geq ij(\ln n)^{2}/n,$ if and only if
$\omega\in(\mathbf{T}_{i,j})^{-1}(F_{n}^{(n)})$ and $ij\leq n/(\ln n)^{2}$.

Then, the probability of $A_{n}$ can be computed as:%
\[
\mathbf{P}(A_{n})=\mathbf{P}(\bigcup\nolimits_{D}\mathbf{T}_{i,j}^{-1}%
(F_{n}^{(n)})=P(\bigcup\nolimits_{D}T_{i,j}^{-1}(F_{n})),
\]
where the union and have indexes in the set $D=\{ij\leq(n-1)/(\ln n)^{2};\ln
n\leq i,j\leq n-1\}.$ By Lemma \ref{Rohling2}, it follows that%
\[
\mathbf{P}(A_{n})=P(F_{n})\sum_{\ln n\leq j\leq n-1}\text{ \ \ }\sum_{\ln
n\leq i\leq(n-1)/j(\ln n)^{2}}1\sim\frac{n\ln n}{2n^{2}(\ln n)^{2}}=\frac
{1}{2n\ln n}.
\]
Therefore%
\[
\sum_{n\geq2}\mathbf{P}(A_{n})=\sum_{n\geq2}\frac{1}{2n\ln n}=\infty.
\]
By the second Borel-Cantelli lemma, $\mathbf{P}(A_{n}$ i.o.$)=1$. This means
that almost all $\omega\in\mathbf{\Omega}$ belong to an infinite number of
$A_{n}$. Whence, taking into account (\ref{basic}), for almost all $\omega
\in\mathbf{\Omega}$ and every positive $B,$%
\begin{equation}
\lim\sup_{i\wedge j\rightarrow\infty}\mathbf{P}^{\omega}(f_{m}\circ
\mathbf{T}_{i,j}/ij\geq B)=1. \label{f untight}%
\end{equation}
Define now%
\begin{equation}
U_{0,0}=\sum_{n\geq2}f_{n}\text{ \ and \ }U_{i,j}=\sum_{n\geq2}f_{n}%
\circ\mathbf{T}_{i,j}. \label{def U}%
\end{equation}
Let us estimate the Luxembourg norm of $U_{0,0}$ in the Orlicz space generated
by the convex function $g(x)=x\ln^{1-\varepsilon}(1+x)$ for $x>0$,
$0<\varepsilon<1$. For each $n\in N$
\[
||f_{n}||_{g}=\inf_{\lambda}\{\lambda:E(\frac{f_{n}}{\lambda}\ln
^{1-\varepsilon}(1+\frac{f_{n}}{\lambda}))\leq1\}\text{ }.
\]
By the definition of $f_{n}$, we have
\begin{align*}
E(\frac{f_{n}}{\lambda}\ln^{1-\varepsilon}(1+\frac{f_{n}}{\lambda}))  &
=P(F_{n})\frac{n}{\lambda\ln^{2}n}\ln^{1-\varepsilon}(1+\frac{n}{\lambda
\ln^{2}n})\\
&  =\frac{1}{2\lambda n\ln^{2}n}\ln^{1-\varepsilon}(1+\frac{n}{\lambda\ln
^{2}n}).
\end{align*}
From this identity we see that, after some computations, that for $n$
sufficiently large%
\[
||f_{n}||_{g}\leq\frac{1}{n\ln^{1+\varepsilon/2}n}.
\]
Clearly, we have%
\begin{equation}
||U_{0,0}||_{g}\leq\sum_{n\geq2}||f_{n}||_{g}<\infty. \label{orlicz norm}%
\end{equation}
It remains to note that, by definition (\ref{def U}), $U_{i,j}\geq f_{n}%
\circ\mathbf{T}_{i,j}$. Therefore, by (\ref{f untight}) we also have for
almost all $\omega\in\mathbf{\Omega}$ and every positive $B,$%
\[
\lim\sup_{i\wedge j\rightarrow\infty}\mathbf{P}^{\omega}(U_{i,j}/ij\geq B)=1
\]
and the conclusion of this lemma follows by letting $B\rightarrow\infty$.
$\square$

\bigskip

\textbf{End of proof of Theorem \ref{Counterex}}

\bigskip

On the space constructed in Lemma \ref{NTight} define the independent random
variables $\xi_{i,j}^{^{\prime}}(\omega_{1},\omega_{2},...)=\xi_{i,j}%
(\omega_{1})~$and the random variables $X_{i,j}=\xi_{i,j}^{\prime}%
U_{i-1,j-1}^{1/2},$ where $(U_{i,j})_{i,j\in Z}$ and $(\xi_{i,j})_{i,j\in Z}$
are as in Lemma \ref{NTight}. Note that $(X_{i,j})_{i,j\in Z}$ is a sequence
of orthomartingale differences with respect to $\mathbf{\ }\prod_{m=1}%
^{\infty}\mathcal{F}_{i,j}^{(m)}$, where $\mathcal{F}_{i,j}^{(m)}$ are
independent copies of $\mathcal{F}_{i,j}$. According to Lemma \ref{NTight} for
$P-$almost all $\omega\in\mathbf{\Omega}$ we have
\[
\lim_{B\rightarrow\infty}\lim\sup_{i\wedge j\rightarrow\infty}P^{\omega
}(|X_{i,j}|/\sqrt{ij}\geq B)=1.
\]
If we assume now that $(S_{n,m}/\sqrt{nm})_{n,m\geq1}$ satisfies the quenched
limit theorem (or it is "quenched" tight), because%
\[
U_{i-1,j-1}^{1/2}=|X_{i,j}|\leq|S_{i,j}|+|S_{i-1,j}|+|S_{i,j-1}|+|S_{i-1,j-1}%
|,
\]
then necessarily the field $(|X_{m,m}|/\sqrt{nm})_{n,m\geq1}$ should be tight
under $P^{\omega},$ for almost all $\omega$, which leads to a contradiction.
Note that, by (\ref{orlicz norm}), for any $0<\varepsilon<1$ we have
$EX_{0,0}^{2}\ln^{1-\varepsilon}(1+|X_{0,0}|)<\infty.$ For this example
$EX_{0,0}^{2}\ln(1+|X_{0,0}|)=\infty,$ since otherwise the quenched result
follows by Theorem \ref{QCLTmart}. $\square$

\section{Quenched functional CLT \label{functional CLT}}

In this section we formulate the functional CLT, which holds under the same
conditions as in Theorem \ref{QCLTmart}. For $(s,t)\in\lbrack0,1]^{2},$ we
introduce the stochastic process
\[
W_{n,v}(t,s)=\frac{1}{\sqrt{nv}}S_{[nt],[vs]}.
\]
We shall establish the following result. Denote by $(W(t,s))_{(t,s)\in
\lbrack0,1]^{2}}$ the standard $2$-dimensional Brownian sheet.

\begin{theorem}
\label{Th funct CLT}Under the setting of Theorem \ref{QCLTmart}, if we assume
that $E(X_{0,0}^{2})<\infty$ then, for $P$-almost all $\omega,$ the sequence
of processes $(W_{n,n}(t,s))_{n\geq1}$ converges in distribution in
$D([0,1]^{2})$ endowed with the uniform topology to $\sigma W(t,s),$ under
$P^{\omega}$. If we assume now that (\ref{moment condition}) holds, then for
$P$-almost all $\omega,$ the sequence $(W_{n,v}(t,s))_{n,v\geq1}$ converges in
distribution to $\sigma W(t,s),$ as $n\wedge v\rightarrow\infty$ under
$P^{\omega}$.
\end{theorem}

\textbf{Proof of Theorem \ref{Th funct CLT}}

\bigskip

Let us first prove the second case, when $n\wedge v\rightarrow\infty.$ As
usual, the proof of this theorem involves two steps, namely the proof of the
convergence of the finite dimensional distributions to the corresponding ones
of the standard $2$-dimensional Brownian sheet and tightness.

For proving tightness we shall verify the moment condition given in relation
(3) in Bickel and Wichura (1971) and then the tightness follows from Theorem 3
in the same paper. To verify it is enough to compute the $4-$th moment of an
increment of the process $W_{n,v}(t,s)$ on the rectangle $A=[t_{1}%
,t_{2})\times\lbrack s_{1},s_{2}).$ That is $E(\Delta^{4}(A))$ where
\[
\Delta(A)=\frac{1}{\sqrt{nv}}\sum\nolimits_{i=[nt_{1}]}^{[nt_{2}]-1}%
\sum\nolimits_{j=[vs_{1}]}^{[vs_{2}]-1}X_{i,j}.
\]
By applying Burkholder's inequality twice consecutively, and taking into
account that the variables are bounded by $C,$ for a positive constant $K$ we
obtain
\[
E^{\omega}(\Delta^{4}(A))\leq KC^{4}(t_{2}-t_{1})^{2}(s_{2}-s_{1})^{2}%
=KC^{4}\mu^{2}(A),
\]
where $\mu$ is the Lebesgue measure on $[0,1]^{2}.$ If $B$ is a neighboring
rectangle of $A$, by the Cauchy-Schwatz inequality we have
\[
E^{\omega}(\Delta^{2}(A)\Delta^{2}(B))\leq KC^{4}\mu(A)\mu(B).
\]
Therefore the moment condition in relation (3) in Bickel and Wichura (1971) is
verified with $\gamma=4$ and $\beta=2.$

The proof of the convergence of finite dimensional distribution follows, up to
a point, the proof of the corresponding result in Cuny et al. (2015), which
will be combined with the method of proof in Theorem \ref{QCLTmart}. As
explained in Subsection 3.2 in Cuny et al. (2015), in order to establish the
convergence of the finite dimensional distributions, we have to show that for
$P-$almost all $\omega\in\Omega,$ and for any partitions $0\leq t_{1}%
\leq...\leq t_{K}\leq1$ and $0\leq s_{1}\leq...\leq s_{K}\leq1,$ we have%

\begin{equation}
\frac{1}{\sqrt{nv}}\sum\nolimits_{k=1}^{K}\sum\nolimits_{\ell=1}^{K}a_{k,\ell
}\sum\nolimits_{i=[nt_{k-1}]}^{[nt_{k}]-1}\sum\nolimits_{j=[vs_{\ell-1}%
]}^{[vs_{\ell}]-1}X_{i,j}\Rightarrow N(0,\Gamma)\text{ under }P^{\omega
},\text{ } \label{fdd}%
\end{equation}
where $\Gamma=\sigma^{2}\sum\nolimits_{k=1}^{K}\sum\nolimits_{\ell=1}%
^{K}a_{k,\ell}^{2}(t_{k}-t_{k-1})(s_{\ell}-s_{\ell-1}).$ Since we have proved
tightness in $C([0,1]^{2}),$ we know that any subsequence contains one which
is converges in distribution to a continuous process. Therefore, without
restricting the generality we can restrict ourselves to partitions with
rational ends which form a countable set.

In order to establish this weak convergence we follow step by step the proof
of Theorem \ref{QCLTmart}. We shall just mention the differences. The first
step is to decompose $X_{i,j}$ as in formula (\ref{truncation}) and to show
the negligibility of the term containing $X_{i,j}^{"}.$ This is the only step
where we need different moment conditions according to whether indexes in the
sum are restricted or not. By using simple algebraic manipulations, the
triangle inequality along with Theorem 3.2\ in Billingsley (1999), we can
easily see that this term is negligible $P$-a.s. for the convergence in
$D([0,1]^{2})$ endowed with the uniform topology, if, for every $\varepsilon
>0$
\[
\lim_{A\rightarrow\infty}\text{ }\lim\sup_{n\wedge v\rightarrow\infty}%
P_{0,0}(\max_{1\leq i\leq n}\max_{1\leq j\leq v}|\sum\nolimits_{k=1}^{i}%
\sum\nolimits_{\ell=1}^{j}\mathcal{P}_{k,\ell}(X_{k,\ell}^{"})|>\varepsilon
\sqrt{nv})=0\text{ a.s.}%
\]
But by using Cairoli's maximal inequality for orthomartinagles (see Theorem
2.3.1 in Khoshnevisan, 2002, p. 19) the proof is reduced to showing
(\ref{negli1}), which was already established in proof of Theorem
\ref{QCLTmart}. Without loss of generality we redenote $\mathcal{P}%
_{i,j}(X_{i,j}^{\prime})$ by $X_{i,j}$ and assume that it is bounded by a
positive constant $C$. We continue the steps of the proof in Theorem
\ref{QCLTmart} and we shall verify the conditions of Theorem \ref{cltGH} with
the exception that we replace $F_{i,v}$ in definition (\ref{defFiv}) by%
\[
F_{k,i,v}=\frac{1}{\sqrt{v}}\sum\nolimits_{\ell=1}^{K}a_{k,\ell}%
\sum\nolimits_{j=[vs_{\ell-1}]}^{[vs_{\ell}]-1}X_{i,j},
\]
where $[nt_{k-1}]\leq i\leq\lbrack nt_{k}]-1;$ $1\leq k\leq K.$ We also
replace $\sigma^{2}$ by $\eta_{k}^{2}=\sigma^{2}\sum\nolimits_{\ell=1}%
^{K}a_{k,\ell}^{2}(s_{\ell}-s_{\ell-1})$ and $h_{m,N_{0}}$ by
\[
h_{k,m,N_{0}}=\sup_{v>N_{0}}E_{0,0}|\frac{1}{m}\sum\nolimits_{i=0}%
^{m-1}F_{k,i,v}^{2}-\eta_{k}^{2}|.
\]
For instance, let us convince ourselves that (\ref{erg2}) holds. Indeed by the
ergodicity of $S$ and the fact that the variables are bounded%
\[
\lim_{v\rightarrow\infty}E_{\infty,0}F_{k,i,v}^{2}=\lim_{v\rightarrow\infty
}\frac{1}{v}E_{\infty,0}(\sum\nolimits_{\ell=1}^{K}a_{k,\ell}\sum
\nolimits_{j=[vs_{\ell-1}]}^{[vs_{\ell}]-1}X_{i,j}^{2})=\eta_{k}^{2}.
\]

After we verify the conditions of Theorem \ref{cltGH} for the triangular array
of martingale differences $(F_{k,i,v})_{[nt_{k-1}]\leq i\leq\lbrack
nt_{k}]-1;\text{ }1\leq k\leq K}$ , we obtain the result in (\ref{fdd}) by
applying the CLT in Theorem \ref{cltGH}. \ $\square$

\section{Quenched functional CLT via coboundary decomposition \label{cob}}

Now we indicate a larger class than the orthomartingale, which satisfies a
quenched functional CLT. A fruitful approach is to approximate $S_{m,n}$ by an
orthomartingale $M_{n,m}$ in a norm that makes possible to transport the
quenched functional CLT given in Theorem \ref{Th funct CLT}. Such an
approximation is of the form: for every $\varepsilon>0,$%
\begin{equation}
\lim\sup_{n\wedge v\rightarrow\infty}P^{\omega}(\max_{1\leq k\leq n,1\leq
\ell\leq v}|S_{k,\ell}-M_{k,\ell}|>\varepsilon\sqrt{nv})=0\text{
}\ \text{a.s.} \label{ortho approx}%
\end{equation}

The random fields we consider can be decomposed into an orthomartingale and a
generalized coboundary and therefore satisfy (\ref{ortho approx}). This type
of orthomartingale approximation, so called martingale-coboundary
decomposition, was introduced for random fields by Gordin (2009) and studied
by El Machkouri and Giraudo (2016), Giraudo (2018) and Voln\'{y} (2018).

\begin{definition}
\label{defcobdec}We say that a random field $(X_{i,j})_{i,j\in Z}$, defined by
(\ref{def X}),\ adapted to the commuting filtration $(\mathcal{F}%
_{i,j})_{i,j\in Z},$ defined by (\ref{def F}),\ admits a martingale-coboundary
decomposition if%
\begin{equation}
X_{0,0}=m_{0,0}+(1-\hat{T})m_{0,0}^{\prime}+(1-\hat{S})m_{0,0}^{"}+(1-\hat
{T})(1-\hat{S})Y_{0,0}, \label{decomp}%
\end{equation}
with $m_{0,0}$ an orthomartingale difference (satisfying (\ref{ortho})),
$m_{0,0}^{\prime}$ a martingale difference in the second coordinate and
$m_{0,0}^{"}$ a martingale difference in the first coordinate. All these
functions are $\mathcal{F}_{0,0}-$measurable.
\end{definition}

We shall obtain the following generalization of Theorem \ref{Th funct CLT}:

\begin{theorem}
\label{CoQCLT}Let us assume that the decomposition (\ref{decomp}) holds with
all the variables square integrable and $S$ $($or $T)\ $is ergodic. Then for
almost all $\omega\in\Omega,$%
\begin{equation}
\frac{1}{n}S_{[nt],[ns]}\Rightarrow|c|W(t,s)\text{ under }P^{\omega}\text{
when }n\rightarrow\infty, \label{clt1}%
\end{equation}
where $(W(t,s))_{(t,s)\in\lbrack0,1]^{2}}$ is the standard $2$-dimensional
Brownian sheet and $c^{2}=E(m_{0,0}^{2})$. If we assume that all the variables
involved in the decomposition (\ref{decomp}) satisfy (\ref{moment condition})
then, for almost all $\omega\in\Omega,$
\begin{equation}
\frac{1}{(nv)^{1/2}}S_{[nt],[vs]}\Rightarrow|c|W(t,s)\text{ under }P^{\omega
}\text{ when }n\wedge v\rightarrow\infty. \label{CLT 2}%
\end{equation}

\end{theorem}

It should be noted that Giraudo (2018) have shown that if
\begin{equation}
\sup_{n,v\geq0}E((E_{0,0}(S_{n,v}))^{2})<\infty, \label{con1}%
\end{equation}
then the decomposition (\ref{decomp}) holds and all the variables are in
$L_{2}.$ As a matter of fact this is also a necessary condition for
(\ref{decomp}). The only condition specific to $L_{2}$ needed for his proof is
the reflexivity of $L_{2}.$ Since the Orlicz space $L_{\varphi}$ generated by
the function
\[
\varphi(x)=x^{2}\log(1+x):[0,\infty)\rightarrow\lbrack0,\infty)
\]
is reflexive (see Theorem 8 in Milnes (1957)), the proof of Theorem 2.1 in
Giraudo is also valid in this context. It follows that if
\begin{equation}
\sup_{n,v\geq0}E(\varphi(|E_{00}(S_{n,v})|))<\infty, \label{con2}%
\end{equation}
then the decomposition in (\ref{decomp}) holds all the functions are in
$L_{\varphi}.$ The reciprocal is also true.

As a matter of fact, by combining Theorem \ref{CoQCLT} with this result we
deduce the following corollary:

\begin{corollary}
\label{corco}Let us assume that the random field $(X_{i,j})_{i,j\in Z}$,
defined by (\ref{def X}),\ adapted to the commuting filtration $(\mathcal{F}%
_{i,j})_{i,j\in Z},$ defined by (\ref{def F}),~satisfies (\ref{con1}). Then
$\lim_{n\wedge v\rightarrow\infty}(nv)^{-1}E(S_{n,v}^{2})=c^{2}$. If in
addition we assume that $S$ $($or $T)\ $is ergodic, then for $P-$almost all
$\omega\in\Omega,$ (\ref{clt1}) holds. Also, if condition (\ref{con2}) is
satisfied, then for $P-$almost all $\omega\in\Omega,$ (\ref{CLT 2}) holds.
\end{corollary}

\bigskip

\textbf{Proof of Theorem \ref{CoQCLT}}

\bigskip

Consider first that the indexes $n$ and $m$ are varying independently. Denote
by $m_{i,j}=m_{0,0}\circ T_{i,j}$ and $M_{k,\ell}=\sum_{i=0}^{k-1}\sum
_{j=0}^{\ell-1}m_{i,j}.$

We shall establish (\ref{ortho approx}). A simple computation shows that
$(S_{k,\ell}-M_{k,\ell})/\sqrt{nv}$ is the sum of the following three terms:%
\[
\frac{1}{\sqrt{nv}}\sum_{i=0}^{k-1}\sum_{j=0}^{\ell-1}\hat{T}^{i}\hat{S}%
^{j}(I-\hat{T})m_{0,0}^{\prime}=\frac{1}{\sqrt{nv}}\sum_{j=0}^{\ell-1}\hat
{S}^{j}(m_{0,0}^{\prime}-\hat{T}^{k}m_{0,0}^{\prime})=R_{1}(k,\ell),
\]%
\[
\frac{1}{\sqrt{nv}}\sum_{i=0}^{k-1}\sum_{j=0}^{\ell-1}\hat{T}^{i}\hat{S}%
^{j}(I-\hat{S})m_{0,0}^{"}=\frac{1}{\sqrt{nv}}\sum_{i=0}^{k-1}\hat{T}%
^{i}(m_{0,0}^{"}-\hat{S}^{\ell}m_{0,0}^{"})=R_{2}(k,\ell),
\]%
\[
\frac{1}{\sqrt{nv}}\sum_{i=0}^{k-1}\sum_{j=0}^{\ell-1}\hat{T}^{i}\hat{S}%
^{j}(I-\hat{T})(I-\hat{S})Y_{0,0}=\frac{1}{\sqrt{nv}}(I-\hat{S}^{\ell}%
)(I-\hat{T}^{k})Y_{0,0}=R_{3}(k,\ell).
\]
In order to treat the last term, note that%
\[
\max_{1\leq k\leq n,1\leq\ell\leq v}|R_{3}(k,\ell)|\leq\frac{4}{\sqrt{nv}%
}\text{ }\max_{0\leq i\leq n}\text{ }\max_{0\leq j\leq v}|Y_{i,j}|.
\]
Let $A$ be a positive integer. By truncation at the level $A$ we obtain the
following bound%
\[
\frac{1}{nv}\max_{0\leq i\leq n}\max_{0\leq j\leq v}|Y_{i,j}|^{2}\leq
\frac{A^{2}}{nv}+\frac{1}{nv}\sum_{i=0}^{n}\sum_{j=0}^{v}Y_{i,j}^{2}%
I(|Y_{i,j}|>A).
\]
Because of the stationarity and the fact that in the second part of Theorem
\ref{CoQCLT}\textbf{ }we imposed condition\textbf{ (}\ref{moment condition}%
\textbf{)}, by the ergodic theorem for stationary random fields (see Theorem
1.1 in Ch.6, Krengel (1985)) it follows that for every $A,$
\[
\lim_{n\wedge v\rightarrow\infty}\frac{1}{nv}\sum_{i=0}^{n}\sum_{j=0}%
^{v}Y_{i,j}^{2}I(|Y_{i,j}|>A)=E(Y_{0,0}^{2}I(|Y_{0,0}|>A)).
\]
Therefore $\lim_{A\rightarrow\infty}\lim_{n\wedge v\rightarrow\infty}%
|R_{3}(n,v)|=0$ $P-$a.s. By Fubini's theorem, it follows that the limit is $0$
also under $P^{\omega}$, for almost all $\omega$.

The terms $R_{1}(k,\ell)$ and $R_{2}(k,\ell)$ are treated similarly, with
small differences. Let us treat the first one only. It is convenient to
truncate at a positive number $A$.$\ $Let%
\begin{align*}
m_{j,k}^{\prime}  &  =m_{j,k}^{\prime}I(|m_{j,k}^{\prime}|\leq A)-E_{j,k-1}%
m_{j,k}^{\prime}I(|m_{j,k}^{\prime}|\leq A)+\\
m_{j,k}^{\prime}I(|m_{j,k}^{\prime}|  &  >A)-E_{j,k-1}m_{j,k}^{\prime
}I(|m_{j,k}^{\prime}|>A).
\end{align*}
We shall use the following bound:%
\begin{gather*}
E_{0,0}\max_{1\leq k\leq n,1\leq\ell\leq v}R_{1}^{2}(k,\ell)\leq2E_{0,0}%
\max_{1\leq k\leq n,1\leq\ell\leq v}(\sum_{j=0}^{\ell-1}m_{j,k}^{\prime}%
)^{2}\leq\\
8A^{2}v+2E_{0,0}\max_{1\leq k\leq n,1\leq\ell\leq v}(\sum_{j=0}^{\ell
-1}m_{j,k}^{\prime}I(|m_{j,k}^{\prime}|>A)-E_{j,k-1}m_{j,k}^{\prime}%
I(|m_{j,k}^{\prime}|>A))^{2}\\
\leq8A^{2}v+2\sum_{k=1}^{n}E_{0,0}\max_{1\leq\ell\leq v}(\sum_{j=0}^{\ell
-1}m_{j,k}^{\prime}I(|m_{j,k}^{\prime}|>A)-E_{j,k-1}m_{j,k}^{\prime}%
I(|m_{j,k}^{\prime}|>A))^{2}.
\end{gather*}
Now, by the Doob's maximal inequality
\begin{gather*}
\frac{1}{nv}E_{0,0}\max_{1\leq k\leq n,1\leq\ell\leq v}R_{1}^{2}(k,\ell)\\
\leq\frac{8A^{2}}{n}+\frac{2}{nv}\sum_{k=1}^{n}\sum_{j=0}^{v-1}E_{0,0}%
(m_{j,k}^{\prime}I(|m_{j,k}^{\prime}|>A)-E_{j,k-1}m_{j,k}^{\prime}%
I(|m_{j,k}^{\prime}|>A))^{2}\\
\leq\frac{8A^{2}}{n}+\frac{4}{nv}\sum_{k=1}^{n}\sum_{j=0}^{v-1}E_{0,0}%
(m_{j,k}^{\prime}I(|m_{j,k}^{\prime}|>A))^{2}\\
=\frac{8A^{2}}{n}+\frac{4}{nv}\sum_{k=1}^{n}\sum_{j=0}^{v-1}Q_{1}^{j}Q_{2}%
^{k}[(m_{0,0}^{\prime})^{2}I(|m_{0,0}^{\prime}|>A)].
\end{gather*}
We let $n\wedge v\rightarrow\infty$\ and we use Theorem 1.1 in Ch. 6 of
Krengel (1985). It follows that, for every $A$
\[
\lim_{n\wedge v\rightarrow\infty}\frac{1}{nv}E_{0,0}\max_{1\leq k\leq
n,1\leq\ell\leq v}R_{1}^{2}(k,\ell)=E(m_{0,0}^{\prime})^{2}I(|m_{0,0}^{\prime
}|>A).
\]
Then, we let $A\rightarrow\infty.$ This completes the proof of
(\ref{ortho approx}). The result follows by using the second part of Theorem
\ref{Th funct CLT} along with Theorem 3.2 in Billingsley (1999). Now for the
situation $n=m\rightarrow\infty,$ the proof is similar with the difference
that we use Theorem 3.5 in Ch. 6 in Krengel (1985)\ instead of Theorem 1.1 in
the same chapter together with the first part of Theorem \ref{Th funct CLT}.
$\square$

\begin{remark}
If we take $Y_{0,0}$, in the martingale-coboundary decomposition
(\ref{decomp}),\ to be the function $U_{0,0}^{1/2}$ found in the proof of
Lemma \ref{NTight}, then for almost all $\omega$,%
\[
R_{3}(n,v)=\frac{1}{\sqrt{nv}}\sum_{i=0}^{n-1}\sum_{j=0}^{v-1}\hat{T}^{i}%
\hat{S}^{j}(I-\hat{T})(I-\hat{S})Y_{0,0}%
\]
does not converge to $0$ in probability $P^{\omega}$ when $n\wedge
v\rightarrow\infty.$ Therefore if only the existence of the second moment is
assumed or even if $EY_{0,0}^{2}\ln^{1-\varepsilon}(1+|Y_{0,0}|)<\infty$ for
some $0<\varepsilon<1,$ this coboundary could spoil the quenched weak
convergence. This is in sharp contrast with the dimension $1$. Recall that in
dimension $1$, when we have a martingale-coboundary decomposition $X_{0}%
=D_{0}+G_{0}-\hat{T}G_{0}$ with $D_{0}$ a martingale difference and $G_{0}\in
L_{2},$ then the coboundary $G_{0}-\hat{T}G_{0}$ does not spoil the quenched
invariance principle (see Theorem 8.1 in Borodin and Ibragimov (1994), which
is due to Gordin and Lifshits). In higher dimension, in general, we need
stronger moment conditions not only for martingale differences but also for
the cobounding function $Y_{0,0}$.
\end{remark}

\bigskip

\section{The case of d-indexed random field \label{section gend}}

In this section we formulate our results and indicate their proofs for random
fields indexed by $Z^{d}$ with $d>2.$ The proofs are based on induction
arguments. When we add on unrestricted $d$-dimensional rectangles the moment
conditions will depend on $d$. By $\mathbf{u=}(u_{1},u_{2},...,u_{d})$ we
denote elements of $Z^{d}$. Let us suppose that $\mathbf{T}=(T_{i})_{1\leq
i\leq d}$ are $d$ commuting, invertible, measure preserving transformations
from $\Omega$ to $\Omega$ and let $\mathcal{F}_{\mathbf{0}}$ be a sub-sigma
field of $\mathcal{K}$. For all $\mathbf{u}\in Z^{d}$ define $\mathcal{F}%
_{\mathbf{u}}=\mathbf{T}^{-\mathbf{u}}(\mathcal{F}_{\mathbf{0}}),$ where
$\mathbf{T}^{-\mathbf{u}}$ is the following composition of operators:
$\mathbf{T}^{-\mathbf{u}}=\prod\nolimits_{i=1}^{n}T_{i}^{-u_{i}}.$ Assume the
filtration is coordinatewise increasing and commuting, in the sense that for
any integrable variable we have $E_{\mathbf{u}}E_{\mathbf{a}}X=E_{\mathbf{a}%
\wedge\mathbf{u}}X,$ where $\mathbf{a}\wedge\mathbf{u}$ means coordinatewice
minimum and we used the notation $E_{\mathbf{u}}X=E(X|\mathcal{F}_{\mathbf{u}%
})$. We introduce the stationary field by starting with a $\mathcal{F}%
_{\mathbf{0}}-$measurable function $X_{\mathbf{0}}:\Omega\rightarrow R$ and
then define the random field $X_{\mathbf{k}}(\omega)=X_{\mathbf{0}}%
(\mathbf{T}^{\mathbf{k}}(\omega))=X_{\mathbf{0}}(T_{1}^{k_{1}}\circ...\circ
T_{d}^{k_{d}}).$ The operator $\mathbf{\hat{T}}$ is defined on $L_{2}$ as
$\mathbf{\hat{T}(}f\mathbf{)=}f\circ\mathbf{T.}$ For the filtration
$(\mathcal{F}_{\mathbf{u}})_{\mathbf{u}\in Z^{d}},$ defined as above,\ we call
the random field $(X_{\mathbf{u}})_{\mathbf{u\in Z}^{d}}$ \ orthomartingale
difference if $E(X_{\mathbf{u}}|\mathcal{F}_{\mathbf{i}})=0$ when at least one
coordinate of $\mathbf{i}$ is strictly smaller that the corresponding
coordinate of $\mathbf{u}$. We also use the notation $\mathbf{i}%
\leq\mathbf{u,}$ where the inequality is coordinatewise and $|\mathbf{n}%
|=n_{1}\cdot...\cdot n_{d}.$ Finally denote $S_{\mathbf{n}}=%
{\displaystyle\sum\nolimits_{\mathbf{0}\leq\mathbf{i}\leq\mathbf{n-1}}}
X_{\mathbf{i}}.$ In this context we have:

\begin{theorem}
\label{trmd}Assume that there is an integer $i$, $1\leq i\leq d$ such that
$T_{i}$ is ergodic and $X_{\mathbf{0}}$ is square integrable, $E(X_{\mathbf{0}%
}^{2})=\sigma^{2}$. Then, for $P-$almost all $\omega\in\Omega,$%
\[
\frac{1}{n^{d/2}}S_{(n,n,...,n)}\Rightarrow\sigma W(t_{1},...,t_{d})\text{
under }P^{\omega}\text{ when }n\rightarrow\infty\text{ }.
\]
In addition, if $E[X_{\mathbf{0}}^{2}\log^{d-1}(1+|X_{\mathbf{0}}|)]<\infty$,
then for almost all $\omega\in\Omega,$
\[
\frac{1}{|\mathbf{n}|^{1/2}}S_{\mathbf{(}n_{1},n_{2},...,n_{d})}%
\Rightarrow\sigma W(t_{1},...,t_{d})\text{ under }P^{\omega}\text{ when }%
\min_{1\leq i\leq d}n_{i}\rightarrow\infty.
\]

\end{theorem}

\begin{remark}
\label{remark d}Both Theorems \ref{Th funct CLT} and \ref{CoQCLT} as well as
Corollary \ref{corco} also hold for the multi-indexed random field
$(X_{\mathbf{u}})_{\mathbf{u\in Z}^{d}}$ defined above.
\end{remark}

We shall indicate how to prove these results by induction. We shall follow
step by step the proof of Theorem \ref{QCLTmart} with the following
differences. Without restricting the generality, let us assume that the
operator $T_{i}$ is ergodic for an $i$, $2\leq i\leq d.$ We define now the
$d$-dimensional projection operators. By using the commutative property of the
filtrations it is convenient to define:%

\[
\mathcal{P}_{\mathbf{u}}(X)=\mathcal{P}_{\mathbf{u}_{1}}\circ\mathcal{P}%
_{\mathbf{u}_{2}}\circ...\circ\mathcal{P}_{\mathbf{u}_{d}}(X),
\]
where%
\[
\mathcal{P}_{\mathbf{u}_{j}}(Y)=E(Y|\mathcal{F}_{\mathbf{u}})-E(Y|\mathcal{F}%
_{\mathbf{u}_{j}}).
\]
Above we used the notation $\mathbf{u}_{j}$ for a vector which has the same
coordinates as $\mathbf{u}$ with the exception of the $j$-th coordinate, which
is $u_{j}-1$. For instance when $d=3,$ $\mathcal{P}_{\mathbf{u}_{2}%
}(Y)=E(Y|\mathcal{F}_{u_{1},u_{2},u_{3}})-E(Y|\mathcal{F}_{u_{1},u_{2}%
-1,u_{3}}).$~We can easily see that, by using the commutativity property of
the filtration, this definition is a generalization of the case $d=2$. We note
that, by using this definition of $\mathcal{P}_{\mathbf{u}}(X),$ the
truncation argument in Theorem \ref{QCLTmart} remains unchanged if we replace
the index set $Z^{2}$ with $Z^{d}.$ We point out the following two differences
in the proof of Theorem \ref{trmd}. One difference is that, for the validity
of the limit in (\ref{erg}) when $\min_{1\leq i\leq d}n_{i}\rightarrow\infty,$
in order to apply the ergodic theorem for Dunford-Schwartz operators, conform
to Ch. 6 Theorem 2.8 and Theorem 1.1 in Krengel (1985), we have to assume that
$E[X_{0,0}^{2}\log^{d-1}(1+|X_{0,0}|)]<\infty$. After we reduce the problem to
the case of bounded random variables, we proceed with the proof of the CLT by
induction. More precisely, we write the sum in the form%
\[
\frac{1}{|\mathbf{n}|^{1/2}}S_{\mathbf{(}n_{1},n_{2},...,n_{d})}=\frac
{1}{n_{1}^{1/2}}\sum_{k_{1}=0}^{n_{1}-1}F_{k_{1},(n_{2},n_{3},...,n_{d}%
)}\text{ ,}%
\]
with%
\[
F_{k_{1},(n_{2},n_{3},...,n_{d})}=\frac{1}{(n_{2}\cdot...\cdot n_{d})^{1/2}%
}\sum_{\mathbf{k\in}B}X_{\mathbf{k}}\text{ ,}%
\]
where the sum is taken on the set $B=\{(0,...,0)\leq(k_{2}...k_{d})\leq
(n_{2}-1,...,n_{d}-1)\}.$ Because one operator is ergodic, according to the
induction hypothesis, $F_{k,(n_{2},n_{3},...,n_{d})}\Rightarrow N(0,\sigma
^{2})$ under $P^{\omega}$ for almost all $\omega,$ and we can replace
(\ref{erg2}) by
\[
\lim\frac{1}{n_{2}...n_{d}}E_{\infty,0,...0}\sum_{B}X_{\mathbf{k}}^{2}%
=\sigma^{2}\text{ a.s. when }\min(n_{2},...,n_{d})\rightarrow\infty.
\]

\section{Examples\label{Appl}}

We shall give examples providing new results for linear and Volterra random
fields with i.i.d. innovations. Let $d$ be an integer $d>1.$ Denote by
$\mathbf{t=(}t_{1},t_{2},...,t_{d})$ and let $W(\mathbf{t})$ be the standard
$d$-dimensional Brownian sheet.

\begin{example}
\label{exlinear} Let $(\xi_{\mathbf{n}})_{\mathbf{n}\in Z^{d}}$ be a random
field of independent, identically distributed random variables, which are
centered and have finite second moment. Let $(a_{\mathbf{n}})_{\mathbf{n}\in
Z^{d}}$ be a sequence of real numbers such that $\sum_{\mathbf{j}%
\geq\mathbf{0}}a_{\mathbf{j}}^{2}<\infty$. Define%
\[
X_{\mathbf{k}}=\sum_{\mathbf{j}\geq\mathbf{0}}a_{\mathbf{j}}\xi_{\mathbf{k}%
-\mathbf{j}}.
\]
Assume that
\begin{equation}
\sup_{\mathbf{n}\geq\mathbf{1}}\sum_{\mathbf{i}\geq\mathbf{0}}b_{\mathbf{n}%
,\mathbf{i}}^{2}<\infty,\text{ where }b_{\mathbf{n,i}}=\sum_{\mathbf{0}%
\leq\mathbf{k\leq n-1}}a_{\mathbf{k}+\mathbf{i}}. \label{lin}%
\end{equation}
Then, if $\mathbf{n}=(n,n,...,n),$ for $P-$almost all $\omega$
\begin{equation}
\frac{1}{n^{d/2}}S_{[(\mathbf{n-1)\cdot t]}}\Rightarrow|c|W(\mathbf{t})\text{
under }P^{\omega}\text{ when }n\rightarrow\infty. \label{ex1}%
\end{equation}
If we assume now that $E(\xi_{\mathbf{0}}^{2}\log^{d-1}(1+|\xi_{\mathbf{0}%
}|))<\infty,$ then for $P-$almost all $\omega$
\begin{equation}
\frac{1}{|\mathbf{n}|^{1/2}}S_{[(\mathbf{n-1)\cdot t]}}\Rightarrow
|c|W(\mathbf{t})\text{ under }P^{\omega}\text{ when }\min(n_{1},...,n_{d}%
)\rightarrow\infty, \label{ex2}%
\end{equation}
where $\mathbf{n}=(n_{1},...,n_{d})$.
\end{example}

\textbf{Proof of Example \ref{exlinear}}.

\bigskip

For this case we take $\mathcal{F}_{\mathbf{n}}=\sigma(\xi_{\mathbf{u}%
},\mathbf{u}\leq\mathbf{n}).$ Let us note first that the variables are square
integrable and well defined. We also have
\[
E(S_{\mathbf{n}}|\mathcal{F}_{\mathbf{0}})=\sum_{\mathbf{0}\leq\mathbf{k\leq
n-1}}\sum_{\mathbf{j}\leq\mathbf{0}}a_{\mathbf{k}-\mathbf{j}}\xi_{\mathbf{j}}%
\]
and therefore%
\[
E(E^{2}(S_{\mathbf{n}}|\mathcal{F}_{\mathbf{0}}))=\sum_{\mathbf{i}%
\geq\mathbf{0}}(\sum_{\mathbf{0}\leq\mathbf{k\leq n-1}}a_{\mathbf{k}%
+\mathbf{i}})^{2}E(\xi_{\mathbf{1}}^{2}).
\]
The result follows for $S_{\mathbf{n}}$ by applying the first part of
Corollary \ref{corco}.

On the other hand, by the Rosenthal inequality for independent random
variables (see relation 21.5 in Burkholder (1973)), applied with the function
$\varphi(x)=x^{2}\log^{d-1}(1+|x|),$ there is a positive constant $C$ such
that
\[
E(\varphi(|E(S_{\mathbf{n}}|\mathcal{F}_{\mathbf{0}})|))\leq C\varphi\left(
\left(  \sum_{\mathbf{i}\geq\mathbf{0}}b_{\mathbf{n,i}}^{2}E(\xi_{\mathbf{1}%
}^{2})\right)  ^{1/2}\right)  +C\sum_{\mathbf{i}\geq\mathbf{0}}E(\varphi
(|b_{\mathbf{n,i}}\xi_{\mathbf{0}}|)),
\]
which is bounded under condition (\ref{lin}). Indeed, condition (\ref{lin})
implies that $\sup_{\mathbf{n}\geq\mathbf{1}}\sup_{\mathbf{i}\geq\mathbf{0}%
}|b_{\mathbf{n},\mathbf{i}}|<\infty,$ and then, after simple algebraic
manipulations we can find a positive constant $K$ such that
\[
E(\varphi(|b_{\mathbf{n,i}}\xi_{\mathbf{0}}|))\leq Kb_{\mathbf{n,i}}%
^{2}(E(\varphi(|\xi_{\mathbf{0}}|))+E(\xi_{\mathbf{0}}^{2})).
\]
It remains to apply the second part from Corollary \ref{corco} and Remark
\ref{remark d} in order to obtain the second part of the example. $\square$

\bigskip

Another class of nonlinear random fields are the Volterra processes, which
play an important role in the nonlinear system theory.

\begin{example}
\label{Volterra}Let $(\xi_{\mathbf{n}})_{\mathbf{n}\in Z^{d}}$ be a random
field of independent random variables, identically distributed, centered and
with finite second moment. Define%
\[
X_{\mathbf{k}}=\sum_{(\mathbf{u},\mathbf{v)}\geq(\mathbf{0},\mathbf{0}%
)}a_{\mathbf{u},\mathbf{v}}\xi_{\mathbf{k-u}}\xi_{\mathbf{k-v}},
\]
where $a_{\mathbf{u},\mathbf{v}}$ are real coefficients with $a_{\mathbf{u}%
,\mathbf{u}}=0$ and $\sum_{\mathbf{u,v}\geq\mathbf{0}}a_{\mathbf{u,v}}%
^{2}<\infty$. Denote%
\[
c_{\mathbf{u},\mathbf{v}}(\mathbf{j})=\sum\limits_{\mathbf{0}\leq\mathbf{k\leq
j-1}}a_{\mathbf{k+u,k+v}}\text{.}%
\]
Assume that%
\begin{equation}
\sup_{\mathbf{j\geq1}}\sum\limits_{\mathbf{u}\geq\mathbf{0,v}\geq
\mathbf{0,u}\neq\mathbf{v}}c_{\mathbf{u},\mathbf{v}}^{2}(\mathbf{j})<\infty.
\label{volt}%
\end{equation}
Then the quenched functional CLT in (\ref{ex1}) holds. If in addition we
assume that $E(\xi_{\mathbf{0}}^{2}\log^{d-1}(1+|\xi_{\mathbf{0}}|))<\infty$,
then the quenched functional CLT in (\ref{ex2}) holds for sums of variables in
a general $d$-dimensional rectangle.
\end{example}

\textbf{Proof of Example \ref{Volterra}}.

\bigskip

For this case we consider the sigma algebras as in Example \ref{exlinear}. We
start from the following estimate
\[
E(S_{\mathbf{j}}|\mathcal{F}_{\mathbf{0}})=\sum_{(\mathbf{u},\mathbf{v)}%
\geq(\mathbf{0},\mathbf{0})}\sum\limits_{\mathbf{0}\leq\mathbf{k\leq j-1}%
}a_{\mathbf{k+u},\mathbf{k+v}}\xi_{-\mathbf{u}}\xi_{-\mathbf{v}}%
=\sum_{(\mathbf{u},\mathbf{v)}\geq(\mathbf{0},\mathbf{0})}c_{\mathbf{u}%
,\mathbf{v}}(\mathbf{j})\xi_{-\mathbf{u}}\xi_{-\mathbf{v}}.
\]
Since by our conditions $c_{\mathbf{u},\mathbf{u}}=0,$ by Tonelli theorem we
obtain
\begin{align*}
E(E^{2}(S_{\mathbf{j}}|\mathcal{F}_{\mathbf{0}}))  &  =\sum\limits_{\mathbf{u}%
\geq\mathbf{0,v}\geq\mathbf{0,u}\neq\mathbf{v}}(c_{\mathbf{u},\mathbf{v}}%
^{2}(\mathbf{j})+c_{\mathbf{u},\mathbf{v}}(\mathbf{j})c_{\mathbf{v}%
,\mathbf{u}}(\mathbf{j}))E(\xi_{\mathbf{u}}\xi_{\mathbf{v}})^{2}\\
&  \leq2\sum\limits_{\mathbf{u}\geq\mathbf{0,v}\geq\mathbf{0,u}\neq\mathbf{v}%
}(c_{\mathbf{u},\mathbf{v}}^{2}(\mathbf{j})+c_{\mathbf{v},\mathbf{u}}%
^{2}(\mathbf{j}))E(\xi_{\mathbf{u}}\xi_{\mathbf{v}})^{2}\\
&  \leq4\sum\limits_{\mathbf{u}\geq\mathbf{0,v}\geq\mathbf{0,u}\neq\mathbf{v}%
}c_{\mathbf{u},\mathbf{v}}^{2}(\mathbf{j})(E(\xi_{\mathbf{0}}^{2}))^{2}.
\end{align*}
The first result of this theorem follows by applying the first part of
Corollary \ref{corco} via Remark \ref{remark d}.

On the other hand, by a moment inequality for $U$-statistics based on the
decoupling procedures, (see Relation 3.1.3. in Gin\'{e} et al., 2000), we
obtain for $\varphi(x)=x^{2}\log^{d-1}(1+|x|),$
\[
E(\varphi(|E(S_{\mathbf{j}}|\mathcal{F}_{\mathbf{0}})|))\leq CE\varphi\left(
\sum_{(\mathbf{u},\mathbf{v)}\geq(\mathbf{0},\mathbf{0})}c_{\mathbf{u}%
,\mathbf{v}}(\mathbf{j})\xi_{-\mathbf{u}}\xi_{-\mathbf{v}}^{\prime}\right)  ,
\]
where $(\xi_{\mathbf{n}}^{\prime})_{\mathbf{n}\in Z^{d}}$ in an independent
copy of $(\xi_{\mathbf{n}})_{\mathbf{n}\in Z^{d}}$ and $C$ is a positive
constant. Now, we apply Rosenthal inequality, given in relation 21.5 in
Burkholder (1973), and we find a constant $C^{\prime}>0$ such that
\begin{align*}
E(\varphi(|E(S_{\mathbf{j}}|\mathcal{F}_{\mathbf{0}})|))  &  \leq C^{\prime
}\varphi\left(  \left(  \sum\limits_{\mathbf{u}\geq\mathbf{0,v}\geq
\mathbf{0,u}\neq\mathbf{v}}c_{\mathbf{u},\mathbf{v}}^{2}(\mathbf{j})E^{2}%
(\xi_{\mathbf{1}}^{2})\right)  ^{1/2}\right) \\
&  +C^{\prime}\sum_{(\mathbf{u},\mathbf{v)}\geq(\mathbf{0},\mathbf{0}%
)}E\varphi\left(  |c_{\mathbf{u},\mathbf{v}}(\mathbf{j})\xi_{-\mathbf{u}}%
\xi_{-\mathbf{v}}^{\prime}|\right)  .
\end{align*}
Note that, by (\ref{volt}), we have $\sup_{\mathbf{u},\mathbf{v\geq0,j\geq1}%
}|c_{\mathbf{u},\mathbf{v}}(\mathbf{j})|<\infty.$ Also, because $\xi
_{-\mathbf{u}}$ and $\xi_{-\mathbf{v}}^{\prime}$ are independent and
identically distributed, by the properties of $\varphi,$ we can find positive
constants such that
\begin{align*}
&  E\varphi\left(  |c_{\mathbf{u},\mathbf{v}}(\mathbf{j})\xi_{-\mathbf{u}}%
\xi_{-\mathbf{v}}^{\prime}|\right) \\
&  \leq Kc_{\mathbf{u},\mathbf{v}}^{2}(\mathbf{j})[E(\varphi(\xi_{\mathbf{0}%
}))E(\xi_{\mathbf{0}}^{2})+(E(\xi_{\mathbf{0}}^{2}))^{2}]\leq K^{\prime
}c_{\mathbf{u},\mathbf{v}}^{2}(\mathbf{j}).
\end{align*}
It remains to note that condition (\ref{volt}) implies condition (\ref{con2})
and then to apply the second part of Corollary \ref{corco} and Remark
\ref{remark d}.

\begin{remark}
In Examples \ref{exlinear} and \ref{Volterra} the innovations are i.i.d.
fields. However, the property (\ref{com}) for the filtration is not restricted
to filtrations generated by independent random variables. For example, we can
take as innovations the random field $(\xi_{n,m})_{n.m\in Z}$ having as
columns independent copies of a stationary and ergodic martingale differences
sequence. In this case the filtration generated $(\xi_{n,m})_{n,m\in Z}$ is
also commuting. As a matter of fact a commuting filtration could be generated
by a stationary random field $(\xi_{n,m})_{n,m\in Z}$ where the columns are
independent, i.e. $\bar{\eta}_{m}=(\xi_{n,m})_{n\in Z}$ are independent.
\end{remark}

\section{Auxiliary results \label{AUX}}

The following is a variant of Theorem 1 in G\"{a}nssler and H\"{a}usler (1979)
(see also G\"{a}nssler and H\"{a}usler, 1986, pages 315--317).

\begin{theorem}
\label{cltGH}Assume that $(D_{n,k})_{1\leq k\leq n}$ is a triangular array of
martingale differences adapted to an increasing filtration $(\mathcal{F}%
_{n,k})_{k}$. Assume that for all $q\ $rational numbers in $[0,1],$%
\begin{equation}
\sum_{k=1}^{[nq]}D_{n,k}^{2}\rightarrow^{P}\sigma^{2}q \label{GH}%
\end{equation}
and $\max_{1\leq k\leq n}|D_{n,k}|$ is uniformly integrable. Then
$S_{[nt]}\Rightarrow\sigma W(t),$ where $S_{[nt]}=\sum_{k=1}^{[nt]}D_{n,k}$
and $W(t)$ is a standard Brownian measure. In particular $S_{n}\Rightarrow
N(0,\sigma^{2})$.
\end{theorem}

As a matter of fact, condition (\ref{GH})\ in Theorem 1 in G\"{a}nssler and
H\"{a}usler (1979) is formulated for all reals $t\in\lbrack0,1].$ We notice
however that if (\ref{GH}) holds for any $q$ rational number in $[0,1]$ then
it also holds for any $t\in\lbrack0,1].$ To see it fix $t,$ $t\in\lbrack0,1]$
and let $q_{1}$ and $q_{2}$ be two rational numbers such that $q_{1}\leq t\leq
q_{2}.$ Then, by using monotonicity, note that
\[
\sum\nolimits_{k=1}^{[nq_{1}]}D_{n,k}^{2}-\sigma^{2}q_{2}\leq\sum
\nolimits_{k=1}^{[nt]}D_{n,k}^{2}-\sigma^{2}t\leq\sum\nolimits_{k=0}%
^{[nq_{2}]}D_{n,k}^{2}-\sigma^{2}q_{1}%
\]
and therefore%
\[
|\sum\nolimits_{k=1}^{[nt]}D_{n,k}^{2}-\sigma^{2}t|\leq\max_{i=1,2}%
|\sum\nolimits_{k=0}^{[nq_{i}]}D_{n,k}^{2}-\sigma^{2}q_{i}|+(q_{2}%
-q_{1})\sigma^{2}.
\]
By using the hypothesis (\ref{GH}), and the fact that the rational numbers are
dense in $R$ it follows that (\ref{GH}) holds for any $t\in\lbrack0,1]$.

\bigskip

We mention now Theorem 34.2 (v) in Billingsley (1995). Further reaching
results including comments of the sharpness of the result below can be found
in Argiris and Rosenblatt (2006).

\begin{theorem}
\label{Bill}Assume that the sequence of random variables $(X_{n})_{n\geq0}$
converges a.s. to $X$ and there is an integrable and positive random variable
$Y$ such that $|X_{n}|\leq Y$ a.s. for all $n\geq0$. Let $\mathcal{F}$ be a
sigma algebra. Then the sequence\ $(E(X_{n}|\mathcal{F}))_{n\geq0}$ converges
a.s. to $E(X|\mathcal{F}).$
\end{theorem}

The following is a result in Katznelson and Weiss (1972) known under the name
of Rokhlin lemma for amenable actions\textbf{.}

\begin{theorem}
\textbf{\label{Rohling}} Let $(\Omega,\mathcal{K},P)$ be a nonatomic
probability space and let $T$ be a measure preserving action of $Z^{2}$ into
$\Omega$
\[
T:Z^{2}\times\Omega\rightarrow\Omega
\]
that is ergodic\textbf{.} Then, for all $\varepsilon>0$ and $n\in N$, there is
a set $B=B(n,\varepsilon)\in\mathcal{K}$ such that for $0\leq i,j\leq n-1$,
$T_{i,j}^{-1}B$ are disjoint for distinct indexes $(i,j)$ and
\[
P(B)\geq\frac{1}{n^{2}}(1-\varepsilon).
\]

\end{theorem}

\textbf{Acknowledgements.} This research was supported in part by the
NSF\ grant DMS--1811373. The first author would like to offer special thanks
to the University of Rouen Normandie, where a large portion of this research
was accomplished during her visit of the Department of Mathematics. The
authors wish to thank the referee for many in-depth comments, suggestions, and
corrections, which have greatly improved the manuscript.

\end{document}